\documentclass[10pt]{article}
\usepackage{amssymb}
\usepackage{amsfonts,amsmath}

\def\Z{\mathbb {Z}}
\def\N{\mathbb {N}}
\def\qed{{\hfill\vrule height 5pt width 5pt depth 0pt}}
\def\phi{\varphi}

\def\Aut{\mathop{\rm Aut}}
\newtheorem{lemma}{Lemma}
\newtheorem{theorem}{Theorem}
\newtheorem{corollary}{Corollary}
\newtheorem{example}{Example}
\begin{document}
\begin{center}
\Large
% TITLE GOES HERE
On Some Generalizations of\\
 Fermat's,
 Lucas's and Wilson's Theorems
\end{center}

\begin{flushright}
Tyler J. Evans \\
Humboldt State University\\
Arcata, CA 95521 USA\\
\verb+evans@humboldt.edu+
\end{flushright}

``Never underestimate a theorem that counts
something!'' -- or so says J. Fraleigh in his classic text
\cite{Fraleigh}.  Indeed, in \cite{Anderson} and \cite{Hausner}, the
authors derive Fermat's (little), Lucas's and
Wilson's theorems, among other results, all
from a single combinatorial lemma. This lemma can be derived by
applying Burnside's theorem to an action by a cyclic group of prime order.
In this note, we
generalize this lemma by applying Burnside's theorem to the corresponding
action by an arbitrary finite cyclic group.  
We revisit the constructions in \cite{Anderson} and
\cite{Hausner} and derive three
divisibility theorems for which the aforementioned classical theorems
are, respectively, the cases
of a prime divisor. Throughout,
$n$ and $p$ denote positive integers with $p$ prime and $\Z_n$ denotes the
cyclic  group of integers
under addition modulo $n$.

\paragraph*{Group Actions and Burnside's Theorem}

By an action of a group $G$ on a set $X$, we mean a homomorphism $G\to\Aut (X)$
where $\Aut(X)$ denotes the group of permutations of $X$.  We
write $gx$ for the image of $x\in X$ under the permutation $X\to X$
induced by $g\in G$.  For each $x\in X$, let $Gx=\{gx\ |\ g\in G\}$ denote
the orbit of $x$ in $X$ and for each $g\in
G$, let $X^g=\{x\in X\ |\ gx=x\}$ denote the set of points fixed by
$g$. If both $G$ and $X$ are finite, 
Burnside's theorem states that the
number of distinct orbits is given by   
\[\frac{1}{|G|}\sum_{g\in G} |X^g|.\]
In particular, $\sum_{g\in G} |X^g|$ is divisible by $|G|$. In the
case that  $G=\Z_n$, for all $g\in\Z_n$, $X^g=X^d$
where $d=(g,n)$ is the greatest common divisor of $g$ and $n$. Each
such $g$ has order $n/d$ and there are $\phi(n/d)$ such elements.   
This observation, together with Burnside's theorem,
gives us the following lemma from which we will derive all of our
results in the sequel.

\begin{lemma} 
\label{mainlemma}
If $X$ is a finite set and $\Z_n\to \Aut(X)$ is a group action, then the number
  of orbits is 
\[\frac{1}{n}\sum_{d|n} \phi\left(\frac{n}{d}\right )|X^d|\]
so that in particular, $\sum_{d|n}\phi(n/d ) |X^d|\equiv 0 \pmod n$.
\qed
\end{lemma}
When $n=p$ is prime, Lemma~\ref{mainlemma} reduces to $|X|\equiv
|X^1|\pmod p$, and this is the combinatorial lemma in
\cite{Anderson} and \cite{Hausner}.

\paragraph*{A Generalization of Fermat's (little) Theorem}

If $a$ is a positive integer and $A=\{1,\dots,a\}$, then 
$\Z_n$ acts on the product
$X=A^n$ by cyclically permuting the coordinates of elements $x\in X$. 
If $g\in\Z_n$ has order $n/d$ then each of the coordinates of $x\in X$
has $n/d$ distinct images under all powers of $g$ so that $g$ fixes
$a^d$ elements of
$X$. Applying Lemma~\ref{mainlemma} gives our first theorem.

\begin{theorem} 
\label{gflt}
For any two positive integers $a$ and $n$,
\begin{equation*}\label{mypoly}
\sum_{d|n}\phi\left( \frac{n}{d}\right ) a^d\equiv 0 \pmod n.
\end{equation*}
\qed
\end{theorem}

\begin{corollary}[Fermat's theorem] For any positive integer $a$, \\
 $a^p\equiv a\pmod p$.
\qed
\end{corollary}
Theorem~\ref{gflt} has appeared numerous times in the literature 
\cite{Fredman,Long,macmahon}.
If $a=1$, then obviously the number of orbits is also equal to 1 and
hence, as a bonus, we recover the well known identity $\sum_{d|n}\phi(d)=n$.

\paragraph*{Wilson's Theorem}

In this section, we revisit an action used in \cite{Anderson}
(in the prime case)
and derive a generalization of Wilson's theorem.  
Let $X$ be the set of all cycles of
length $n$ in the symmetric group $\Aut(\{1,\dots,n\})$. Then
$|X|=(n-1)!$ and the action of
$\Z_n$ on $X$ is defined by 
\[g(a_1,\dots,a_n)=(a_1+g,\dots, a_n+g),\]
where the addition in each position is done modulo $n$. Let $d$ be a
divisor of $n$, $g\in\Z_n$ be an element of order $n/d$,
and let $0,a_2, \dots , a_d\in\Z_n$ be a complete set of
representatives for the set of cosets
$\Z_n / \langle d\rangle$. Define a cycle $\pi=\pi(g,a_2,\dots,a_d)\in X$ by
\begin{equation}
\label{theform}
\pi=(0,a_2, \dots, a_d,g,a_2+g,\dots, a_d+g,
\dots,((n/d)-1)g, \dots, a_d+((n/d)-1)g)
\end{equation}  
where the multiplication is done modulo $n$. There are $\phi(n/d)$
choices for $g$, $(n/d)^{d-1}$ choices for the elements $a_2,\dots,
a_d$ and $(d-1)!$ ways to order them so that the number of cycles of
the form (\ref{theform}) is given by
\[\phi\left(\frac{n}{d}\right ) \left (\frac{n}{d}\right )^{d-1}
(d-1)!.\]
\begin{example}\rm
Let $n=12$, $d=4$, $g=8$, $a_2=9$, $a_3=6$ and $a_4=3$. Then the cycle $\pi$
defined above is
\[\pi = (0,9,6,3,8,5,2,11,4,1,10,7).\]
The reader can verify that $\pi$ is a fixed point only for the
elements in the subgroup $\langle 4\rangle$ of $\Z_{12}$. The cycle
$\pi$ is $1$ of $324=2\cdot 3^3\cdot 3!$ $12$-cycles fixed by the
elements in the subgroup $\langle 4\rangle$. 
\qed
\end{example}
Let $\pi=\pi(g,a_2,\dots, a_d)$. Since $g$ has order $n/d$, we have
$d=kg$ for some 
$1\le k\le (n/d)-1$. It is then easy to see that $d\pi$ is obtained from $\pi$
by cyclically permuting the entries in each position $kd$ spaces to the
left, hence $\pi\in X^d$.
On the other hand, if
$\pi=(a_1,\dots,a_{n})\in X^d$ where $a_1=0$ and $a_{k+1}=d$, 
then $d\pi$ is obtained from $\pi$ by
cyclically permuting the entries $k$ spaces to the left.  It
follows that $k\in\Z_n$ has order $n/d$ so that $d=uk$ for some
$u\in\Z_n$. Therefore
{\it subtracting} $d$ from each entry in $\pi$ a total of $u$ times is
equivalent to moving each entry
{\it right} $d$ spaces.  Since $a_1=0$, this implies $a_{jd+1}=jud$
for all $j=0,\dots (n/d)-1$.  Therefore the order of $ud$ is $n/d$ so that
$a_2\not\in \langle ud\rangle$ and, by similar reasoning,
$a_{jd+2}=a_2+jud$ for all $j=0,\dots (n/d)-1$ exhausting the coset 
$a_2+\langle ud\rangle$.
Continuing, we see that $a_2,\dots, a_d$ represent distinct
cosets in $\Z_n/\langle d \rangle$ and $\pi=\pi(ud,a_2,\dots, a_d)$
has the form (\ref{theform}). We have shown
$|X^g|=\phi(n/d)(n/d)^{d-1}(d-1)!$ so that an application of
Lemma~\ref{mainlemma} gives our second divisibility theorem.
\begin{theorem}
\label{gwt}
For  $n\ge 1$, 
\[\sum_{d|n} \left [\phi\left(\frac{n}{d}\right)\right ]^2\left (\frac{n}{d}\right )^{d-1}
(d-1)! \equiv 0 \pmod n.\]
\qed
\end{theorem}
\begin{corollary}[Wilson's theorem] $(p-1)!\equiv -1\pmod p$.
\qed
\end{corollary}

\paragraph*{Lucas's Theorem}

In this final section, we reanalyze an action used in
\cite{Hausner} (in the prime case) to derive a generalization of
Lucas's  theorem (see Corollary~\ref{lucas} below). 
Let $m,r\ge 0$ and use the division algorithm to write $m=Mn+m_0$ and 
$r=Rn+r_0$ with $0\le m_0,r_0<n$.  For $1\le
k\le n$, let 
\[A_k=\{(k,1), (k,2),\dots, (k,M)\}\ 
\mbox{\rm and let}\ 
B=\{(0,1), (0,2),\dots, (0,m_0)\}.\]
Let $A=A_1\cup A_2\cup\cdots\cup A_n\cup B$ so that $|A|=Mn+m_0=m$. 
Given $C\subseteq A$, let
$C_j=C\cap A_j$ for $1\le j\le n$ and $C_0=C\cap B$ so that 
$C=C_1\cup C_2\cup\cdots\cup
C_n\cup C_0$. If $X$ is the
collection of all $C\subseteq A$ with $|C|=r$, then
$|X|=\binom{m}{r}$. (Note: $\binom{m}{r}=0$ if $m<r$.)
Define $f:A\to A$ by 
\begin{eqnarray*}
f(k,x) & = & (k+1,x) \ \mbox{\rm if $1\le k\le n-1$};\\
f(n,x) & = & (1,x);\\
f(0,x) & = & (0,x),
\end{eqnarray*}
and note easily that $f\in\Aut(A)$. Clearly $f^n$
is the identity map so that the map
$1\mapsto f$ gives an action $\Z_n\to
\Aut(X)$. Moreover, an element $C\in X$ is fixed by $g\in\Z_n$ of
order $n/d$ if and
only if for all $1\le k\le d$, 
$\pi_2(C_k)=\pi_2(C_{lg+k})$ for $l=0,\dots, (n/d)-1$ where $\pi_2$ is
projection onto the second coordinate.  Therefore
\[Rn+r_0=r=|C|=\frac{n}{d}\sum_{k=1}^d|C_k| + |C_0|.\]
But, $0\le r_0,|C_0|<n$, and it follows that there exists
$j\in\{-(d-1),\dots,d-1\}$ such that
\begin{equation}
\label{size}
R=\frac{1}{d}\sum_{k=1}^d |C_k|+\frac{j}{d}\hspace{.5cm}\mbox{\rm
  and}\hspace{.5cm} |C_0|-r_0=(n/d)j.
\end{equation}
Conversely, for all $j\in\{-(d-1),\dots,d-1\}$ and all
 choices of $\alpha_k=|C_k|$ ($1\le k\le d$) that satisfy
(\ref{size}), we can independently choose subsets $C_k\subset A_k$ and
$C_0\subset B$ with $|C_0|=r_0+(n/d)j$, and a unique fixed point of
 $X$ is determined.
If we define the length $||\alpha||_d$ of an element
$\alpha=(\alpha_1,\dots,\alpha_d)\in\N^d$ by
\[||\alpha||_d=\frac{1}{d}\sum_{j=1}^d \alpha_j,\]
then we have shown if $g\in\Z_n$ has order
$n/d$, then
\[|X^g|=\sum_{j=-(d-1)}^{d-1}
  \sum_{||\alpha||_d=\atop R-(j/d)}\binom{M}{\alpha_1}\cdots
\binom{M}{\alpha_d}\binom{m_0}{r_0+(n/d)j}\equiv
0\pmod n\]
Applying Lemma~\ref{mainlemma}, we have our third divisibility theorem.
\begin{theorem} 
\label{glt}
For $n\ge 1$, $m=Mn+m_0$, $r=Rn+r_0$, $0\le m_0,r_0<n$
\[\sum_{d|n}\phi\left (\frac{n}{d}\right )\sum_{j=-(d-1)}^{d-1}
  \sum_{||\alpha||_d=\atop R-(j/d)}\binom{M}{\alpha_1}\cdots
\binom{M}{\alpha_d}\binom{m_0}{r_0+(n/d)j}\equiv
0\pmod n.\]
\qed
\end{theorem} 
\begin{corollary}[Lucas's theorem]\label{lucas} Suppose
\begin{eqnarray*}
m& =& m_kp^k+\cdots+m_1p+m_0;\\
r& =& r_kp^k+\cdots+r_1p+r_0
\end{eqnarray*}
with $0\le m_j,r_j < p$.  Then 
\[\binom{m}{r}\equiv \binom{m_k}{r_k}\cdots\binom{m_1}{r_1}
\binom{m_0}{r_0}\pmod p.\]
\end{corollary}
\noindent {\bf Proof.} We will show that if $m=Mp+m_0$, $r=Rp+r_0$,
$0\le m_0,r_0<p$, then 
\[\binom{m}{r}\equiv \binom{M}{R}\binom{m_0}{r_0}\pmod p,\]
leaving the induction for the reader.
Taking $n=p$, Theorem~\ref{glt} gives
\begin{eqnarray*}
\lefteqn{(p-1)\binom{M}{R}\binom{m_0}{r_0}+}\\
& & \sum_{j=-(p-1)}^{p-1}
\sum_{||\alpha||_p=\atop R-(j/p)}\binom{M}{\alpha_1}\cdots
\binom{M}{\alpha_p}\binom{m_0}{r_0+(n/p)j}\equiv 0\pmod p.
\end{eqnarray*}
Selecting subsets $C_k\subseteq A_k$ with $|C_k|=\alpha_k$ ($1\le k\le
p$) and $C_0\subset B$ with $|C_0|=r_0+(n/p)j$ uniquely determines a
subset $C\in X$ provided $||\alpha||_p=R-(j/p)$.  Therefore
\[\sum_{j=-(p-1)}^{p-1}\sum_{||\alpha||_p=\atop R-(j/p)} 
\binom{M}{\alpha_1}\cdots
\binom{M}{\alpha_p}\binom{m_0}{r_0+(n/p)j}=\binom{m}{r}\]
and the proof is complete.
\qed

\paragraph*{Acknowledgment} The author gratefully acknowledges Mr.~Jon Miller
at Salisbury University for his valuable assistance.

\bibliography{references}

\begin{thebibliography}{1}

\bibitem{Anderson}
Peter~G. Anderson, Arthur~T. Benjamin, and Jeremy~A. Rouse.
\newblock Combinatorial proofs of {F}ermat's, {L}ucas's and {W}ilson's
  theorems.
\newblock {\em Amer. Math. Monthly}, 112(3):266--268, March 2005.

\bibitem{Fraleigh}
John~B. Fraleigh.
\newblock {\em A First Course in Abstract Algebra}.
\newblock Addison-Wesley, 4th edition, 1989.

\bibitem{Fredman}
M.L. Fredman, L.E. Mattics, and L.~Carlitz.
\newblock Elementary {P}roblems and {S}olutions: E2242.
\newblock {\em Amer. Math. Monthly}, 78(5):545--546, May 1971.

\bibitem{Hausner}
Melvin Hausner.
\newblock Applications of a simple counting technique.
\newblock {\em Amer. Math. Monthly}, 90(2):127--129, February 1983.

\bibitem{Long}
Calvin~T. Long.
\newblock Problems and {S}olutions: 6468.
\newblock {\em Amer. Math. Monthly}, 92(10):742, December 1985.

\bibitem{macmahon}
P.A. MacMahon.
\newblock Applications of the theory of permutations in circular procession to
  the theory of numbers.
\newblock {\em Proc. London Math. Soc.}, 23:305--313, 1891-2.

\end{thebibliography}

\end{document}